\documentclass[12pt]{article}

\usepackage{url}

\def\A{$\cal{A}$}
\def\B{$\cal{K}$}
\def\C{$\cal{C}$}
\def\Cp{$\cal{C}'$}

\def\cirk{\,{\raisebox{.3ex}{\tiny $\circ$}}\,}

\def\mj{\mathbf{1}}

\begin{document}

\title{\LARGE On Sets of Premises}
\author{{\sc Kosta Do\v sen}
\\[1ex]
{\small Faculty of Philosophy, University of Belgrade, and}\\[-.5ex]
{\small Mathematical Institute, Serbian Academy of Sciences and Arts}\\[-.5ex]
{\small Knez Mihailova 36, p.f.\ 367, 11001 Belgrade,
Serbia}\\[-.5ex]
{\small email: kosta@mi.sanu.ac.rs}}
\date{}
\maketitle

\begin{abstract}
\noindent Conceiving of premises as collected into sets or multisets, instead of sequences, may lead to triviality for classical and intuitionistic logic in general proof theory, where we investigate identity of deductions. Any two deductions with the same premises and the same conclusions become equal. In terms of categorial proof theory, this is a consequence of a simple fact concerning adjunction with a full and faithful functor applied to the adjunction between the diagonal functor and the product biendofunctor, which corresponds to the conjunction connective.
\end{abstract}

\vspace{2ex}


\vspace{1ex}


\vspace{3ex}

\noindent {\small \emph{Acknowledgements.} Work on this paper was
supported by the Ministry of Education, Science and Technological Development of Serbia,
while the Alexander von Humboldt Foundation has, within the frame of Humboldt-Kollegs, supported the
presentation of a talk partly related to the paper at the conference \emph{Proof} in Bern in September 2013. I am grateful to the organizers of
that conference for their hospitality, and to Milo\v s Ad\v zi\' c and Zoran Petri\' c for discussing the ideas of the paper, reading a draft of it and making useful comments.}

\vspace{3ex}

\section{Introduction}

General proof theory addresses the question ``What is a proof?'', or rather ``What is a deduction?''---a deduction being a hypothetical proof, i.e.\ a proof with hypotheses---by dealing with questions related to normal forms for deductions, and in particular with the question of identity criteria for deductions. It deals with the structure of deductions, as exhibited, for example, with the help of the typed lambda calculus in the Curry-Howard correspondence, and not with their strength measured by ordinals, which is what one finds in proof theory that arose out of Hilbert's programme.

Much of general proof theory is the field of categorial proof theory. Fundamental notions of category theory like the notion of adjoint functor, and very important structures like cartesian closed categories, came to be of central concern for logic in that field. Through results of categorial proof theory called coherence results, which provide a model theory for equality of deductions, logic finds new ties with geometry, topology and algebra (see the books \cite{D99}, \cite{DP04} and \cite{DP07}, the more recent introductory survey \cite{D14}, and references therein).

In general proof theory, and in particular in categorial proof theory, one looks for an algebra of deductions, and for that, one concentrates on the operations of this algebra, which come with the inference rules. As an equational theory, the algebra of deductions involves the question of identity criteria for deductions, the central question of general proof theory. (This question may be found, at least implicitly, in Hilbert's 24th problem; see \cite{D15}.)

In categorial proof theory one usually studies a freely generated category of a certain kind equationally presented. This freely generated category is constructed out of syntactical material, as in universal algebra one constructs a freely generated algebra of a certain kind equationally presented by factoring through an equivalence relation on terms. In categories we have partial algebras---the arrow terms out of which the equivalence classes are built have types, their sources and targets---but there is no significant mathematical difference in the construction when compared with what one has in universal algebra (see \cite{DP04}, Chapter~2, in particular in Section 2.3). The objects of this freely constructed categories are propositions, i.e.\ formulae, and the arrows, i.e.\ the equivalence classes of arrow terms, are deductions, i.e.\ equivalence classes of particular derivations, whose sources are premises and whose targets are conclusions. For deductions we have the partial operation of composition and identity deductions (this is essential for them; see \cite{D16}). The categories in question are interesting if they are not preorders, i.e., not all arrows with the same source and the same target are equal. Otherwise, the proof theory is trivial: any deductions with the same premises and the same conclusions become equal.

In terms of categorial proof theory, assuming that premises are collected into sets leads to assuming that for every proposition $A$ we have that $A$ and $A\wedge A$ are isomorphic, where $\wedge$ is the conjunction connective. Isomorphism is understood here as in category theory: there are arrows, i.e.\ deductions, from $A$ to $A\wedge A$ and back, which composed give identity arrows, i.e.\ identity deductions (see \cite{DP12} and references therein). We are led to assume moreover that the associativity and commutativity of conjunction give isomorphisms.

It is not difficult to establish that a category \B\ with binary product $\times$ is a preorder (i.e., any two arrows with the same source and the same target are equal) iff for every object $B$ of \B\ the diagonal arrow from $B$ to $B\times B$ is an isomorphism. It is also not difficult to establish the related fact that \B\ is a preorder iff for every object $B$ of \B\ the first-projection arrow from $B\times B$ to $B$ is equal to the second-projection arrow with the same source and the same target. We will put these facts within a more general categorial context involving adjunction, which should shed light upon them. This is the main goal of this paper, and achieving this goal, together with some related matters (like those in the last section), is the novelty it should bring.

In categorial proof theory the binary connective of conjunction, both classical and intuitionistic, is identified with binary product. So what we mentioned above indicates that assuming that premises are collected into sets leads to preordered categories where objects are propositions, i.e., formulae, and arrows are deductions. This makes the proof theory trivial.

Assuming that premises are collected into multisets leads to assuming that for every proposition $A$ the permutation deduction from $A\wedge A$ to $A\wedge A$, which permutes the two occurrences of $A$ in $A\wedge A$, is equal to the identity deduction from $A\wedge A$ to $A\wedge A$. Although that assumption leads to the same absolute trivialization of the proof theory of classical and intuitionistic logic, the trivialization obtained with it for linear logic is less absolute. The same less absolute trivialization is obtained for the proof theory of relevant logic by assuming moreover that $A$ and $A\wedge A$ are isomorphic. (What is this less absolute trivialization will be explained in the last section, at the end of the paper.)

\section{Sequents} Gentzen's sequents are expressions of the form $A_1,\ldots ,A_n\vdash B_1,\ldots ,B_m$ where $A_1,\ldots ,A_n,B_1,\ldots ,B_m$ are formulae of an object language, like a language of propositional logic or a first-order language. Instead of the turnstile $\vdash$ Gentzen writes $\rightarrow$ (which is more commonly used nowadays for the binary connective of implication; we use it below, as usual, for separating the sources and targets of arrows in categories), for $A$ and $B$ he uses Gothic letters, and for $n$ and $m$ Greek letters (see \cite{G35}, Section I.2.3). The natural numbers $n$ and $m$ may also be zero; when $n$ is zero $A_1,\ldots ,A_n$ is the empty word, and analogously for $m$ and $B_1,\ldots ,B_m$. For what we have to say in this paper we could restrict ourselves to sequents where $m$ is one.

The comma in sequents is an auxiliary symbol that serves to separate formulae in sequences, which however is not essential. A sequent could as well be $A_1\ldots A_n\vdash B_1\ldots B_m$, but it could be difficult, though not impossible, to see where $A_i$ ends and $A_{i+1}$ begins in the sequence $A_iA_{i+1}$. Instead of $p,p,p\wedge q\vdash r,r$ we would have the less perspicuous $ppp\wedge q\vdash rr$, which however is not ambiguous. (It becomes more perspicuous when we do not omit the outermost parentheses of formulae, as in $pp(p\wedge q)\vdash rr$.)

So Gentzen's sequents may be conceived as expressions of the form $\Gamma\vdash\Delta$ where $\Gamma$ and $\Delta$ are finite, possibly empty, sequences of formulae. (Capital Greek letters as schemata in sequents originate from Gentzen; see op. cit.) To economize upon considerations involving structural rules, a number of authors seem to think they are improving upon Gentzen if they take in ${\Gamma\vdash\Delta}$ that $\Gamma$ and $\Delta$ are not \emph{sequences} of formulae, but other sorts of finite collections of formulae, \emph{multisets} (i.e.\ sets of occurrences) of formulae or \emph{sets} of formulae. One should note immediately that with that $\Gamma\vdash\Delta$ seizes to be a word of a formal language, as usually conceived. If $\Gamma$ and $\Delta$ are multisets or sets, then $\Gamma\vdash\Delta$ is not a \emph{sequence} of symbols. It could be conceived as a triple $(\Gamma, \vdash, \Delta)$, in which case $\vdash$ is not essential. A sequent could be identified with the ordered pair $(\Gamma, \Delta)$.

Such a move is not without its dangers. We are not usually interested in particular sequents, but in sequent schemata. In other words, we want our sequents to be closed under substitution. If we continue using the notation $A_1,\ldots ,A_n\vdash B_1,\ldots ,B_m$ as it is usually done, and if $A_1,\ldots ,A_n$ and $B_1,\ldots ,B_m$ are conceived as sets of formulae, then by substituting $p$ for $q$ from the sequent $p,q\vdash p$ we obtain as a substitution instance $p\vdash p$, and the application of the structural rule of thinning on the left is transformed into the figure on the right:
\[
\frac{p\vdash p}{p,q\vdash p} \hspace{10em} \frac{p\vdash p}{p\vdash p}
\]
which does not look like an application of thinning.

This may be embarrassing, but need not be calamitous. It becomes really dangerous when we are interested not only in provability, but in proofs, i.e.\ not only in deducibility, but in deductions, and try to characterize identity of deductions. If we do that in category theory, i.e.\ in categorial proof theory, making the move that corresponds to switching from sequences to multisets or sets is dangerous, and may result in collapse. It may trivialize matters: any two deductions with the same premises and the same conclusions will be equal. The corresponding categories will be preorders.

\section{Sequents in categorial proof theory}
In categorial proof theory sequents of the simple kind $A\vdash B$, where $A$ and $B$ are single formulae, give the types of arrows $f\!:A\rightarrow B$, with $A$ being the source and $B$ the target. In the presence of conjunction $\wedge$, which serves to replace the comma on the left-hand side of Gentzen's sequents, and the constant $\top$ which replaces the empty sequence, together with disjunction $\vee$ and the constant $\bot$ for the right-hand side, we can mimic Gentzen's sequents, as Gentzen himself envisaged in \cite{G35} (Section I.2.4).

Conjunction, classical or intuitionistic (they are the same), corresponds in this perspective to binary product, and $\top$ corresponds to a terminal object (which may be conceived as nullary product). The associativity natural isomorphism for binary product justifies Gentzen's switch to sequences. The matter with the commutativity, i.e.\ symmetry, natural isomorphism for binary product is more tricky. (We deal with it in the last section.) We may identify $A\wedge B$ and $B\wedge A$, but we should nevertheless distinguish the identity deduction for $A\wedge A$ to $A\wedge A$ from the deduction that permutes the two occurrences of $A$.

To the principle of contraction there does not however correspond an isomorphism: $A$ and $A\wedge A$ are not isomorphic, as the sets $A$ and $A\times A$ are not isomorphic if $A$ is finite with more than one member. To assume that contraction is tied to an isomorphism leads to triviality. Matters are analogous with binary coproduct and an initial object, which correspond respectively $\vee$ and $\bot$, but we will concentrate on conjunction and product in this paper.

By conjunction we understand here a binary connective to which a meet operation of a semilattice corresponds algebraically. In the terminology of substructural logics, this is an additive conjunction, for which we have the two natural-deduction rules of  conjunction elimination, the first-projection and second-projection rules that correspond to Gentzen's structural rule of thinning. (We will consider on another occasion what happens if we assume that contraction is an isomorphism for a multiplicative conjunction that does not involve thinning, like a conjunction we find in relevant logic.)

\section{Adjunction with a full and faithful functor}
We will appeal to a result concerning adjoint functors dual to Theorem~1 of Section IV.3 of \cite{ML}. (The proof of this result below will be more direct and simpler than the proof in loc.\ cit.; it will not appeal to the Yoneda lemma.)

Let $F$ be a functor from a category \B\ to a category \A, and $G$ a functor from \A\ to \B, such that $G$ is right-adjoint to $F$, with members of the counit of the adjunction being $\varphi_A\!:FGA\rightarrow A$, for $A$ an object of the category \A, and members of the unit of the adjunction being $\gamma_B\!:B\rightarrow GFB$, for $B$ an object of the category \B\ (we follow the notational conventions of \cite{D99}, Chapter~4, save that we write \B\ instead of $\cal{B}$, to distinguish this category clearly from its object $B$). Then we have the following.

\vspace{2ex}

\noindent {\sc Proposition 4.1.}\quad {\it The functor $F$ is faithful iff for every object $B$ of \B\ the arrow $\gamma_B$ is monic.}

\vspace{2ex}

\noindent {\it Proof.}\quad From left to right we have:
\begin{tabbing}
\hspace{2em}$\gamma_B\cirk g_1=\gamma_B\cirk g_2$ \= $\Rightarrow F\gamma_B\cirk Fg_1=F\gamma_B\cirk Fg_2$, since $F$ is a functor\\*
\>$\Rightarrow \varphi_{FB}\cirk F\gamma_B\cirk Fg_1=\varphi_{FB}\cirk F\gamma_B\cirk Fg_2$\\
\>$\Rightarrow Fg_1=Fg_2$, by a triangular equation of adjunction\\
\>$\Rightarrow g_1=g_2$, since $F$ is faithful.
\end{tabbing}

From right to left we have:
\begin{tabbing}
\hspace{2em}$Fg_1=Fg_2$ \= $\Rightarrow GFg_1\cirk\gamma_C=GFg_2\cirk\gamma_C$\\*
\>$\Rightarrow \gamma_B\cirk g_1=\gamma_B\cirk g_2$, by the naturality of $\gamma$\\
\>$\Rightarrow g_1=g_2$, since $\gamma_B$ is monic. {\it q.e.d.}
\end{tabbing}

For $\mj_B\!:B\rightarrow B$ being an identity arrow we have the following.

\vspace{2ex}

\noindent {\sc Proposition 4.2.}\quad {\it The functor $F$ is full iff for every object $B$ of \B\ there is an arrow $h_B\!:GFB\rightarrow B$ of \B\ such that $\gamma_B\cirk h_B=\mj_{GFB}$.}

\vspace{2ex}

\noindent {\it Proof.}\quad From left to right, the fullness of $F$ implies that for every object $B$ of \B\ there is an arrow $h_B\!:GFB\rightarrow B$ such that $Fh_B=\varphi_{FB}$. We have:
\begin{tabbing}
\hspace{2em}$\gamma_B\cirk h_B$ \=$=GFh_B\cirk \gamma_{GFB}$, by the naturality of $\gamma$\\*
\>$=G\varphi_{FB}\cirk \gamma_{GFB}$, since $Fh_B=\varphi_{FB}$\\
\>$=\mj_{GFB}$, by a triangular equation of adjunction.
\end{tabbing}

From right to left, take $f\!:FB_1\rightarrow FB_2$ and $j=h_{B_2}\cirk Gf\cirk \gamma_{B_1}\!:B_1\rightarrow B_2$. We have:
\begin{tabbing}
\hspace{2em}$Fj$ \= $=Fh_{B_2}\cirk FGf\cirk F\gamma_{B_1}$, since $F$ is a functor\\*
\>$=\varphi_{FB_2}\cirk F\gamma_{B_2}\cirk Fh_{B_2}\cirk FGf\cirk F\gamma_{B_1}$, by a triangular equation of\\
\`adjunction\\
\>$=\varphi_{FB_2}\cirk FGf\cirk F\gamma_{B_1}$, since $\gamma_{B_2}\cirk h_{B_2}=\mj_{GFB_2}$ and $F$ is a functor\\
\>$=f\cirk\varphi_{FB_1}\cirk F\gamma_{B_1}$, by the naturality of $\varphi$\\
\>$=f$, by a triangular equation of adjunction. {\it q.e.d.}
\end{tabbing}

\noindent {\sc Corollary 4.3.}\quad {\it The functor $F$ is full and faithful iff for every object $B$ of \B\ the arrow $\gamma_B$ is an isomorphism.}

\vspace{2ex}

\noindent {\it Proof.}\quad From left to right, we have $h_B$ by the fullness of $F$ and Proposition 4.2. Then we have:
\begin{tabbing}
\hspace{2em}$F(h_B\cirk\gamma_B)$ \=$=Fh_B\cirk F\gamma_B$, since $F$ is a functor\\*
\>$=\varphi_{FB}\cirk F\gamma_B$, since $Fh_B=\varphi_{FB}$\\
\>$=F\mj_B$,
\end{tabbing}
by a triangular equation of adjunction and the functoriality of $F$.
By the faithfulness of $F$ we obtain that $h_B\cirk\gamma_B=\mj_B$.

From right to left, we use Propositions 4.1 and 4.2 from right to left.
\vspace{-2ex}
\begin{tabbing}
\`{\it q.e.d.}
\end{tabbing}

As another corollary we have that if $F$ is full and faithful, then the adjunction is trivial, in the sense that every two canonical arrows of adjunction of the same type, i.e.\ with the same source and the same target, are equal (for trivial adjunctions see \cite{D99}, Sections 4.6.2 and 4.11).

\section{Adjunction with the diagonal functor}
Consider now the diagonal functor $D$ from a category \B\ to the product category $\cal{K}\times \cal{K}$. This functor is always faithful. We have on the other hand the following.

\vspace{2ex}

\noindent {\sc Proposition 5.1.}\quad {\it The functor $D$ is full iff \B\ is a preorder.}

\vspace{2ex}

\noindent {\it Proof.}\quad From left to right, take the arrows $g,g'\!:B_1\rightarrow B_2$ of \B. By the fullness of $D$, for the arrow $(g,g')\!:(B_1,B_1)\rightarrow (B_2,B_2)$ of $\cal{K}\times \cal{K}$ we have an arrow $h$ of \B\ such that $(g,g')=Dh=(h,h)$. Hence $g=g'$.

From right to left, we have for every arrow $(g,g')\!:(B_1,B_1)\rightarrow (B_2,B_2)$ of $\cal{K}\times \cal{K}$ that $g=g'$. Hence $(g,g')=(g,g)=Dg$. {\it q.e.d.}

\vspace{2ex}

The category \B\ has a product biendofunctor $\times$ iff this functor is right adjoint to the diagonal functor $D$. The members $\gamma_B\!:B\rightarrow B\times B$ of the unit of this adjunction, which can be called diagonal arrows, correspond to the contraction arrows of conjunction $w_B\!:B\vdash B\wedge B$. If these arrows are isomorphisms, then, by Proposition 4.2 (and Corollary 4.3 from right to left), the diagonal functor $D$ is full, and, by Proposition 5.1 from left to right, \B\ is a preorder. In proof-theoretical terms, any two deductions with the same premises and the same conclusions are equal. So if we want our proof theory to be nontrivial, we should not permit the contraction arrows $w_B$ to be isomorphisms.

Conversely, if \B\ is a preorder, by Proposition 5.1 from right to left, the diagonal functor is full, and, by Corollary 4.3 from left to right, the arrows $\gamma_B\!:B\rightarrow B\times B$ are isomorphisms. So we have established the following for categories \B\ with binary product $\times$.

\vspace{2ex}

\noindent {\sc Proposition 5.2.}\quad {\it The category \B\ is a preorder iff for every object $B$ of \B\ the diagonal arrow $\gamma_B\!:B\rightarrow B\times B$ is an isomorphism.}

\vspace{2ex}

Note that it is possible in a category \B\ with binary product $\times$ to have for every object $B$ that $B$ and $B\times B$ are isomorphic without the diagonal map being an isomorphism. We will consider that matter in the next section.

We have established Proposition 5.2 without mentioning explicitly the projection arrows, but they are in the background. It is easy to see that \B\ is a preorder iff for every object $B$ of \B\ we have that the first-projection arrow $k^1_{B,B}\!:B\times B\rightarrow B$ is equal to the second-projection arrow $k^2_{B,B}\!:B\times B\rightarrow B$. From left to right this is trivial, while for the other direction we have that if $k^1_{B,B}=k^2_{B,B}$, then for any two arrows $g,g'\!:C\rightarrow B$ of \B\ we have for $\langle g,g'\rangle\!:C\rightarrow B\times B$ that $k^1_{B,B}\cirk \langle g,g'\rangle=k^2_{B,B}\cirk \langle g,g'\rangle$, and hence $g=g'$.

If the diagonal functor $D$ is full, then for the members $(k^1_{B,B}, k^2_{B,B})\!:{(B\times B,B\times B) \rightarrow(B,B)}$ of the counit of the adjunction of $D$ with the biendofunctor $\times$ we have that $k^1_{B,B}=k^2_{B,B}$. So if the projection arrows $k^1_{B,B}$ and $k^2_{B,B}$ are not equal, then $D$ is not full, and, by Proposition 5.1, the category \B\ is not a preorder. That \B\ is not a preorder follows of course immediately from the inequality of $k^1_{B,B}$ and $k^2_{B,B}$, but it is worth seeing how Proposition 5.1 and with it the wider perspective of the preceding section are involved.

This shows how the non-triviality of \B\ hinges on distinguishing the two projection arrows from $B\times B$ to $B$, i.e.\ the two deductions from $B\wedge B$ to $B$ based on the two natural-deduction rules of conjunction elimination (see the end of Section 3).

In a category with binary product we derive $k^1_{B,B}=k^2_{B,B}$ from the assumption that the diagonal arrow $w_B\!:B\rightarrow B\times B$, called $\gamma_B$ above, is an isomorphism in the following manner. We have in any such category that $k^1_{B,B}\cirk w_B= k^2_{B,B}\cirk w_B=\mj_B$, and then, by the isomorphism of $w_B$, we obtain $k^1_{B,B}=k^2_{B,B}$.

When $k^1_{B,B}=k^2_{B,B}$, then $k^1_{B,B}$, i.e.\ $k^2_{B,B}$, is the arrow inverse to $w_B\!:{B\rightarrow B\times B}$, which makes of the diagonal arrow $w_B$ an isomorphism. We have $k^1_{B,B}\cirk w_B=\mj_B$ anyway. We also have:
\begin{tabbing}
\hspace{2em}$w_B\cirk k^1_{B,B}$ \= $=\langle k^1_{B,B},k^1_{B,B}\rangle$, by the naturality of $w_B$\\
\>$=\langle k^1_{B,B},k^2_{B,B}\rangle$, since $k^1_{B,B}=k^2_{B,B}$\\
\>$=\mj_{B\times B}$.
\end{tabbing}

The reduction to triviality brought by assuming that the diagonal arrows $w_B$ are isomorphisms could have been shown by appealing only to these comments about $k^1_{B,B}=k^2_{B,B}$. We preferred however to put the matter in a wider perspective given by the preceding section. This shows that the matters we deal with are not peculiar to the adjunction with the diagonal functor.

The triviality of the adjunctions involving the quantifiers that are considered in \cite{DP09} (Section 1.4) has the same roots involving sets. The universally quantified formula $\forall x A$ interpreted over a domain with two objects $a$ and $b$ stands either for the conjunction $A^x_a\wedge A^x_b$ or the conjunction $A^x_b\wedge A^x_a$. Moreover, $A^x_a\wedge A^x_a$ and $A^x_a$ are not distinguished. Quantifiers involve sets of objects, and neither multisets nor sequences. So from the point of view of general proof theory the assertion that quantifiers are like conjunctions and disjunctions should be taken with a grain of salt.

\section{Collapse with sets and multisets of premises}
We have noted after Proposition 5.2 that it is possible in a category with binary product $\times$ to have for every object $B$ that $B$ and $B\times B$ are isomorphic without the diagonal arrow $w_B\!:B\rightarrow B\times B$ being an isomorphism. The isomorphism may be produced by something else.
To see that, take the category \C\ with binary product freely generated out of a nonempty set $P$ of generating objects, which are usually taken to be syntactical objects, letters, like propositional variables. (The construction of such equationally presented syntactical freely generated categories is described in detail in \cite{DP04}, Chapter~2, in particular in Section 2.3.) Out of \C\ we can build another category with binary product, which we call \Cp, by replacing the objects of \C, which are like propositional conjunctive formulae, by sets of letters, free generators, occurring in these formulae. If $P$ is the singleton $\{p\}$, then we have in \Cp\ a single object $\{p\}$ (i.e.\ $P$ itself). The category \Cp\ is not a preorder if we keep in it the structure of arrows of \C. (That this is possible is shown by the one-object category that is a skeleton of the category of denumerably infinite sets with functions as arrows.\footnote{I am grateful to Zoran Petri\' c for suggesting this example.}) So the diagonal arrow $w_{\{p\}}\!:\{p\}\rightarrow \{p\}$, though its source and target are both $\{p\}$, need not be equal to the identity arrow $\mj_{\{p\}}\!:\{p\}\rightarrow \{p\}$. Analogously, the projection arrows $k^1_{\{p\},\{p\}}\!:\{p\}\rightarrow \{p\}$ and $k^2_{\{p\},\{p\}}\!:{\{p\}\rightarrow \{p\}}$ need not be equal to $\mj_{\{p\}}$, and they need not be mutually equal.

Once we have passed to \Cp\ and understand this category proof-theoretical\-ly, the temptation is however big to take the diagonal arrow $w_{\{p\}}$, which corresponds to contraction, to be the identity arrow $\mj_{\{p\}}$. In proof theory the switch to sets of premises is usually made in order not to be obliged to keep an account of structural rules. If we have to continue keeping this account, what is the advantage of the switch? We are moreover required to keep an account of contraction when it is invisible, when it consists in passing from ${\{p\}}$ to ${\{p\}}$.

When $w_{\{p\}}$ is identified with $\mj_{\{p\}}$ we have collapse, as we have seen in this paper (and in particular in Proposition 5.2). Collapse also ensues when $k^1_{\{p\},\{p\}}=k^2_{\{p\},\{p\}}=\mj_{\{p\}}\!:\{p\}\rightarrow \{p\}$, which seems tempting and natural to assume when the structural rule of thinning is invisible, when it adds a premise we already have, and premises are collected into sets (see the trivial thinning figure mentioned in Section~2).

Suppose now \C\ is a freely generated category with a symmetric monoidal multiplication biendofunctor $\otimes$ (see \cite{ML}, Section VII.7; these categories are called symmetric associative in \cite {DP04}, Chapter~5), which corresponds to the multiplicative conjunction of linear logic, and is tied to collecting premises into multisets. In \C\ we have a symmetry natural isomorphism $c_{B_1, B_2}\!:{B_1\otimes B_2\rightarrow B_2\otimes B_1}$. Now \Cp\ is obtained from \C\ by replacing the objects of \C, i.e.\ propositional formulae, by sets of occurrences of letters in these formulae, i.e.\ multisets of letters occurring in these formulae. As before, \Cp\ is not a preorder if we keep in it the structure of arrows of \C, and in particular continue distinguishing the symmetry isomorphisms $c_{p, p}\!:p\otimes p\rightarrow p\otimes p$ from the identity arrows $\mj_{p\otimes p}\!:p\otimes p\rightarrow p\otimes p$. The corresponding arrows in \Cp\ go from the multiset $\{p,p\}$ to the multiset $\{p,p\}$.

The assumption $c_{B,B}=\mj_{B\otimes B}$ is however dangerous. If $\otimes$ happens to be binary product $\times$, this assumption leads to collapse, as did the isomorphism of the diagonal arrows and the equality of the projection arrows. This is because we have $k^1_{B,B}\cirk c_{B,B}=k^2_{B,B}$ in categories with binary product, where $c_{B,B}$ is defined as $\langle k^2_{B,B},k^1_{B,B}\rangle$. So collecting premises into multisets in classical and intuitionistic logic is as dangerous as collecting them into sets.

The situation is somewhat different in linear logic. We have the assumption $c_{B,B}=\mj_{B\otimes B}$ in a category called $\mathbf{S}'$ in \cite{DP04} (Section 6.5; this category is obtained out of the symmetric associative category freely generated by a nonempty set of generating objects), which is a preorder, and this shows that the assumption $c_{B,B}=\mj_{B\otimes B}$ makes that all the canonical arrows tied to multiplicative conjunction in linear logic that have the same source and the same target are equal. This does not mean however that the assumption will lead in general to preorder.

If among the free generators we do not have only objects, but also different arrows with the same source and the same target, they will not become equal because of $c_{B,B}=\mj_{B\otimes B}$. This is shown by an inductive argument establishing that for every equation derivable for such a freely generated category, for every generating arrow $f$, we cannot have $f$ on one side of the equation without having it on the other.\footnote{I am grateful to Zoran Petri\' c for suggesting this simple argument.} The same holds for relevant logic and categories related to it (see \cite{P02} and \cite{DP07a}).

In these categories we have, as members of a natural transformation corresponding to contraction, diagonal arrows from $B$ to $B\otimes B$, with $\otimes$ corresponding to multiplicative conjunction as above. In addition to $c_{B,B}=\mj_{B\otimes B}$, we are led to assume that these diagonal arrows are isomorphisms if we assume that premises are collected into sets, and we will show in another place how this makes equal all the canonical arrows tied to multiplicative conjunction that have the same source and the same target. This assumption of isomorphism will however not lead in general to preorder, for reasons analogous to those in the preceding paragraph. So, in linear and relevant logic, with premises collected into multisets or sets the danger is still there, but not as big as in classical and intuitionistic logic.

\end{document}